# A kind of dynamic model of prime numbers

Wang Liang[1*], Huang Yan[2]

[1,2](Department of Control Science and Control Engineering, Huazhong University of Science and Technology, WuHan, 430074 P.R.China)

[2](Department of Math, Huazhong University of Science and Technology, WuHan 430074 P.R.China)

**[Abstract]** A dynamic sieve method is designed according to the basic sieve method. It mainly refers to the symbolic dynamics theory. By this method, we could connect the prime system with familiar 'Logistic Mapping'. An interesting discovery is that the pattern of primes could be depicted by a series of orbits of this mapping. Some heuristic proofs for open problems like twin primes are obtained through this relation. This research gives a new viewpoint for the distribution of prime numbers.

**[Keywords]** Prime number, Symbolic dynamics, Logistic Mapping

## 1 Introduction

The detailed distribution of prime number is still a mystery, but some interesting progresses have been obtained in these years.

The distribution of the prime numbers among the integers seems somewhat random, but it's surely produced by simple sieve method. This makes people believe the primes system is a chaos system. So some new methods in chaos and statistic theory are applied to study the pattern of prime numbers. For example in [1][2] the fractal character of primes was investigated, while in [3] the appropriately defined Lyapunov exponents for the distribution of primes were calculated numerically. Power-law behavior in the distribution of primes and correlations in prime numbers were found in [4]. The differences of consecutive prime numbers are also well studied. Besides the well known 3 period oscillations, the same special periodic behavior also appears in Dirichlet classification and increment (difference of difference) of consecutive prime [5][6]. Most of these researches all show interesting numerical feature in prime numbers, but can't deduce the strict theory results.

The theory research for this topic also obtains some new results. By the methods from analytic number theory and ergodic theory, Ben Green and Terence Tao [7][8]showed that for any positive integer k, there exist infinitely many arithmetic progressions of length k consisting only of prime numbers.

Here we prefer the chaos explanations of prime distribution. In this paper, we design a novel dynamic model for prime numbers by symbolic dynamics theory. The symbolic dynamics theory has been successfully applied in some famous problems like 'Smale horseshoe mapping' and 'Three body problem'. Symbolic dynamic is originally a pure math theory, but its application in chaos analyzing makes it become a simple and powerful tool for dynamic system research, especially for discrete system [9]. This theory also allows us to connect some complex systems with simple mappings. This is the reason that we use this theory to build the dynamic model of primes.

Obviously, the periodic behavior is one of most important features of dynamic system. The periodic expressions for prime numbers are described in the following section. To depict the sieve method, we define a new composition rule of symbolic dynamics in the third section. In the fourth section, we connect this dynamic model with Logistic mappings. Interesting heuristic proofs for some open problems are also discussed in this part. The Last section gives a short summary.

## 2 Periodic expressions for prime numbers

---

[*] Corresponding author: guoypm@hust.edu.cn

Eratosthenes Sieve method of two thousand years ago is still the main method to create the prime numbers. Given a list of the numbers between 1 and N, starting with 2, erase all multiples of 2 up to N, other than 2 itself. Call the remaining set $P_2$. Then return to the beginning and taking the first number greater than 2 and erase all of its multiples up to N, again other than the number itself. Repeating this operation, we could get $P_3, P_4, P_5, \cdots P_R$.

If $R = N^{1/2}$, the set $P_R(N)$ will all be prime numbers. We could give a new representation for this method.

Here the nature number is regarded as the points in a line. The value of nature number is expressed by their position in this line. So 'erased' point is marked by one symbol and 'saving' point by other symbol. So we could use only two symbols to describe the sieve method. Here we select 'R,L'. In fact, any two symbols like '1,0' will also do. We first define the prime number by symbols and then describe the sieve method.

We select $LLLLLLLL\cdots\cdots = (L)^n$ to represent the all natural number points $(0,1,2,3,\cdots\cdots)$ in a line.

For prime number 2, all 'L' in the places that could be divided by 2 are 'erased' and replaced by 'R'. Here the first 'L' presenting 0 is also 'erased'. We get:

$RLRLRL\cdots\cdots = (RL)^n$.

For 3, all 'L' in the places that could be divided by 3 are replaced by 'R', We get:

$RLLRLLRLL\cdots\cdots = (RLL)^n$

Thorough this method, prime number could be expressed by a periodic oscillation symbol sequence. This operation is shown in Fig 1.

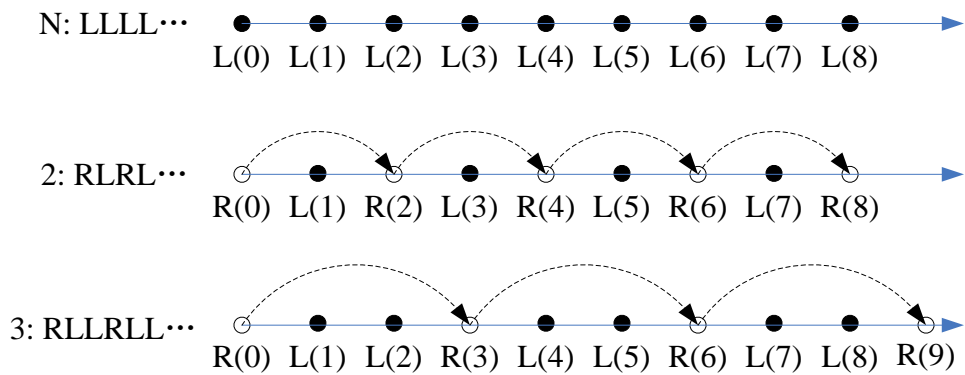

Fig 1. Period expression for prime number

By this method, for any prime number $p_i$, we could get $(RL^{pi})^n$

In this paper, we use the basic symbolic sequence of one period to present the whole period sequence. So the prime number could be described by prime number length symbols:

$2 : M_2 = RL;$

$3 : M_3 = RLL;$

$5 : M_5 = RLLLL;$  $p_i$ is the ith prime number.

$\vdots$

$p_i : M_{Pi} = RL^{pi-1}.$

## 3 Composition rule

After defining the symbolic expression for prime numbers, we could define their composition rule.

For example, we cut all the numbers that could be divided by 2 and 3. Using the method constructing the primes, all 'L' in the place that could be divided by 2 and 3 is 'erased' and replaced by 'R'. We could get $RLRRRL$, it's a 6 period sequence. According to the Sieve rule, the composition of $M_2, M_3$ is very similar with dot product of vectors:

$$M_2 \bullet M_3 = RL \bullet RLL \Leftrightarrow (RL)^n \bullet (RLL)^n = (RLRLRL)^n \bullet (RLLRLL)^n = (RLRLRL \bullet RLLRLL)^n$$

$$\Leftrightarrow (RLRLRL \bullet RLLRLL) = (R \bullet R)(L \bullet L)(R \bullet L)(L \bullet R)(R \bullet L)(L \bullet L)$$

The R presents the 'erased' function. so $(R \bullet R)$ means 'erased' twice: $R \bullet R = R$. $(L \bullet R), (R \bullet L)$ means 'erased' by 2 or 3: $L \bullet R = R \bullet L = R$. Obviously, $L \bullet L = L$ is the 'saving' point.

So $M_2 \bullet M_3 = (R \bullet R)(L \bullet L)(R \bullet L)(L \bullet R)(R \bullet L)(L \bullet L) = RLRRRL$. This operation is shown in Fig 2.

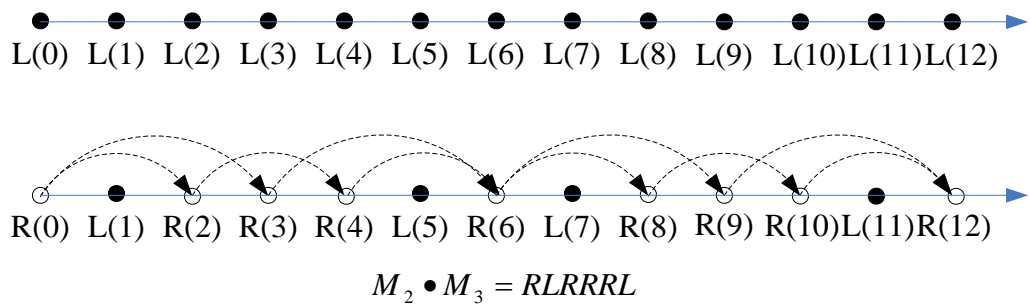

Fig 2 Composition of $M_2, M_3$

There has been a '*' composition rule or called 'DGP' rule in symbolic dynamics theory. Here we need define a new composition rule to describe the sieve method. It's detailed as follows.

If two symbols have the equal length:

$$A = a_1 a_2 a_3 \cdots a_n, B = b_1 b_2 b_3 \cdots b_n$$
$$A \bullet B = (a_1 a_2 a_3 \cdots a_n) \bullet (b_1 b_2 b_3 \cdots b_n) = (a_1 \bullet b_1)(a_2 \bullet b_2)(a_3 \bullet b_3) \cdots (a_n \bullet b_n)$$
$$here, a, b \in \{R, L\}, a \bullet b = \begin{cases} L, a = b = L \\ R, a = b = R \\ R, a \neq b \end{cases}$$

If two symbols are different length, we just need extend their basic period sequence to same length:

$$A = a_1 a_2 a_3 \cdots a_n, B = b_1 b_2 b_3 \cdots b_m$$
$$A \bullet B = (a_1 a_2 a_3 \cdots a_n)^m \bullet (b_1 b_2 b_3 \cdots b_m)^n$$

For any two prime $p_i$, $p_j$, We could prove the period of $M_{Pi} \bullet M_{pj}$ is $p_i \times p_j$ by basic number theory.

Because $M_2 \bullet M_3 = RLRRRL$, according to the sieve method, the place of second 'L' is just the third prime number 5. So $M_5 = RLLLL = f(M_2 \bullet M_3)$, $f(D)$ function selects all the symbols of D before the second 'L' and then turns all the 'R' after the first 'L' to 'L'.

Obviously, for any prime number: $M_{pi+1} = f(M_2 \bullet M_3 \bullet M_5 \bullet \cdots \cdots M_{pi})$, $p_i$ is ith prime number.

Here we mark: $D_1 = M_2, D_2 = M_2 \bullet M_3, \cdots, D_i = M_2 \bullet M_3 \bullet \cdots \bullet M_{pi}$

So 'Eratosthenes sieve method' could be expressed by symbolic method:

$$D_{n+1} = D_n \bullet [f(D_n)], D_1 = RL$$
$$M_{pn} = f(D_n)$$

The advantage of this sieve method is obvious. It turns the division operation of Eratosthenes Sieve into the simple symbol operation. Many methods in symbolic dynamics could also be applied in this research. And most important thing is that we could connect primes system with some easy chaos system through symbols expression. In the following part, we just discuss the relation between primes system and one dimension unimodal mapping.

## 4 Connections with unimodal mapping

Here we are mainly concerned with dynamic feature of prime number system, but not how to produce prime numbers.

We could use $D_i, D_1 = M_2, D_2 = M_2 \bullet M_3, \cdots, D_i = M_2 \bullet M_3 \bullet \cdots \bullet M_{pi}$ to describe the dynamic character of prime number. For example, $D_2 = M_2 \bullet M_3 = RLRRRL$ is 6 period. It contains the period oscillation 'information' of prime number 2 and 3. The dynamical character of whole prime system could be depicted as:

$$D_\infty = M_2 \bullet M_3 \bullet \cdots \bullet M_{pi} \bullet \cdots \cdots$$

Its period $T(D_\infty) = p_1 \times p_2 \times p_3 \times \cdots \cdots \to \infty$. It seems that prime system will transform from period system into a chaos system. The $D_i, i = 1 \to \infty$ could describe this process.

The relation between symbolic dynamics and one dimension unimodal mapping has been well investigated. So we connect the primes system with some simple unimodal mappings to study the pattern of primes. Here we select Logistic mapping $x_{n+1} = 1 - ux_n^2$. The relation of its parameter u and related symbols is shown in Table 1, which could be found in many books for 'Applied symbolic dynamics' [9].

| u | 1.0 | 1.3107 | 1.3815 | 1.40115 | 1.4304 | 1.5437 | 1.754 | 2.0 |
|---|---|---|---|---|---|---|---|---|
| Symbols | $RC$ | $RLRC$ | $RLR^3LRC$ | $R^{*\infty}$ | $RLR^2(RL)^\infty$ | $RLR^\infty$ | $RLC$ | $RL^\infty$ |

Table 1. The relation of parameter u and related sequence

According to the period expression of prime number, $M_2, M_3, \cdots, M_{pi}, \cdots$ are all kneading sequence, which means we could find their corresponding actual orbit from mapping $x_{n+1} = 1 - ux_n^2$. We could also prove their '•' composition $D_1 < D_2 < D_3 < \cdots < D_i < \cdots < D_\infty$ are kneading sequence. So prime system's special road to the chaos could be described by a series of orbits of $x_{n+1} = 1 - ux_n^2$.

For $D_2 = M_2 \bullet M_3 = RLRRRL$, we could get $u(D_2) \approx 1.476$. It's very difficult to get the expression of $D_\infty$ from '•' composition rule. But by our 'sieve method', we know all the nature number will be 'erased' except 1, so $D_\infty = RLRRR\cdots = RLR^\infty$. From Table 1, we get $u(D_\infty) = u(RLR^\infty) \approx 1.5437$. The series orbits $x_{n+1} = 1 - ux_n, 1.476 < u < 1.5437$ could depict the primes' special road to the chaos, which are shown in the bifurcation figure of Logistic mapping [Fig 3].

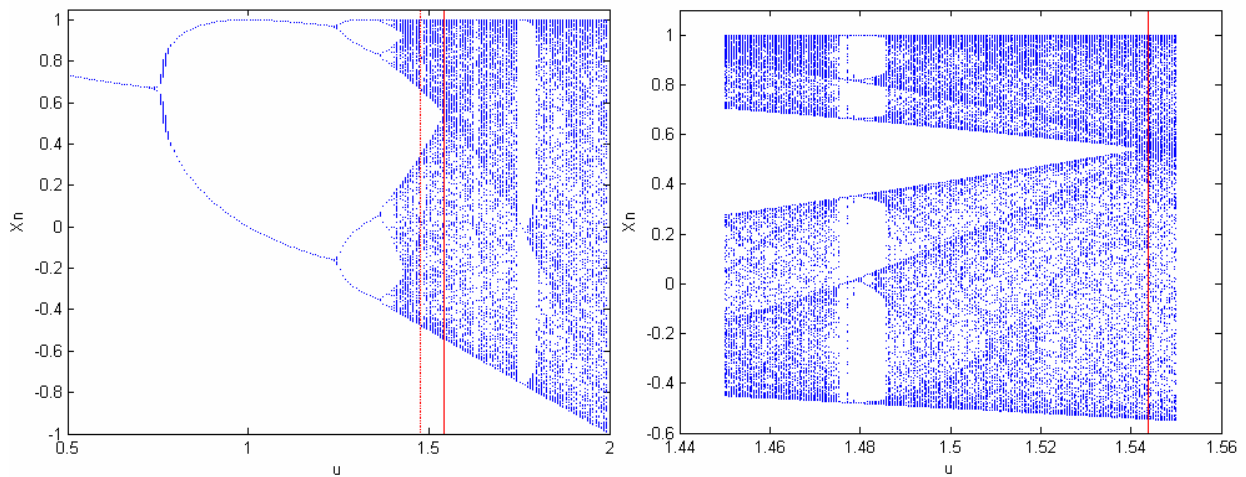

Fig 3 (a), (b) orbits of $x_{n+1} = 1 - ux_n, 1.476 < u < 1.5437$ in bifurcation figure

In Fig 3 (a), the first red line corresponds to $u = 1.476$, the second line is $u = 1.5437$. Fig 3(b) is the detailed figure near $1.476 < u < 1.5437$. The red line corresponds to $u = 1.5437$.

The mapping $x_{n+1} = 1 - u(D_\infty)x_n^2$ could represent the dynamical character of whole prime system. According to the kneading theory, the topological entropy of $D_\infty = RLR^\infty$ is $\ln(2)/2$. The Lyapunov exponent of $x_{n+1} = 1 - u(D_\infty)x_n^2$ is about 0.3406.

Here the prime numbers distribution problem is converted into the research of the orbits near $x_{n+1} = 1 - u(D_\infty)x_n^2$. Some interesting heuristic proof for open problems in number theory could be obtained through this relation:

**(1) Twin primes conjecture**

The $u(D_\infty)$ in the bifurcation figure is the converging point of 'two bands' and 'single band'. So according to the basic principle of ergodic theory, in enough long time, the chaotic orbits near $u(D_\infty)$ could contain arbitrarily amount of symbols segment in the form of d:

$$d = (\underbrace{LRRRRR \cdots R}_{2n})L$$

, here n is natural number.

On the other hand, the periodic orbits of $D_i, i \to \infty$ also lie near the $u(D_\infty)$. Here $\sum u(D_i)$ denotes the symbols expression of chaotic orbits near $u(D_\infty)$. We could construct $\sum u(D_i)$ that 'immensely' approach to $D_i$, which means $\sum u(D_i)$ and $D_i$ have enough long equal head part.

When constructing $D_i = M_2 \bullet M_3 \bullet \cdots \bullet M_{pi}$, we 'erased' all the 'L' point that could be divided by $p_1, p_2, p_3, \cdots p_i$. Through this operation we could get the prime numbers between $(p_i, p_i^2)$. So the first $n, n \approx \dfrac{p_i^2}{2\ln(p_i)} - \dfrac{p_i}{\ln(p_i)} = \dfrac{p_i(p_i - 2)}{2\ln(p_i)}$ number of symbols 'L' in $D_i$ all represent the prime number except the 'L(1)'. We could also say when $i \to \infty, n \to \infty$, the 'L' points in $\sum u(D_i)$ except 'L(1)' are all prime numbers in enough long time.

$\sum u(D_i)$ could contain arbitrarily amount of segments in the form of d in enough long time. That is to say arbitrarily amount of even number less than a certain value will appear in the difference of consecutive prime numbers between $(p_i, p_i^2)$, here $p_i$ is a large prime number. Part of this result could also be described as there are arbitrarily number of twin primes between $(p_i, p_i^2)$. It provides a method to construct the infinity amount of twin primes.

Here we are still not very clear about the approaching method of $D_i, i \to \infty$ and relation between period orbits of $D_i$ and chaos orbits nearby. We will discuss it in the further paper.

**(2) Goldbach conjecture**

Here r(n) denotes the number of representations of an even number n as the sum of two primes. In our paper [2], the computational results suggest the periodic oscillations behavior of series $r(n), n \to \infty$ could be described by the symbolic sequence:

$$M_r = M_3 * M_5 * M_7 * \cdots * M_{pi} * \cdots$$

In this paper, we find the periodic behavior of primes system is:

$$D_P = D_\infty = M_2 \bullet M_3 \bullet \cdots \bullet M_{pi} \bullet \cdots\cdots$$

If defining a special sum operation '$\oplus$' for symbols, we could describe Goldbach Conjecture as:

$$M_r = D_P \oplus D_P$$

The relation of two combination rules '*', '$\bullet$' and '$\oplus$' may crack the relation of multiplication and addition of primes. Now there is still no an efficient method to study the Goldbach conjecture. This symbolic dynamics method may open a possible avenue for success in proving this conjecture.

The Lyapunov exponent of $x_{n+1} = 1 - u(M_r)x_n^2$ is $L(r) \in (0.1983, 0.2119)$. It seems that the primes system is more 'chaotic' than the statistics of sum of two primes. In fact, we could easily find the distribution of primes is almost random but the "Goldbach Comet" is well-regulated.

**(3) Riemann hypothesis**

Riemann zeta function also plays the important role in symbolic dynamics. It connects the period symbol sequence with the chaotic sequence. Some interesting research direction for Riemann hypothesis may be found by analyzing the orbits near $x_{n+1} = 1 - u(D_\infty)x_n^2$.

## 5 Conclusions

This paper builds a novel dynamic model for primes system by symbolic dynamics. We find the pattern of primes system could be depicted by some orbits of Logistic Mapping. Many other kinds of chaos system could also be connected with the primes system through symbolic dynamics. Maybe more interesting results could be obtained. The discussion in this paper is not very precise. Only wish it could spark some new ideas in this area.

"Mathematicians have tried in vain to this day to discover some order in the sequence of prime numbers, and we have reason to believe that it is a mystery into which the mind will never penetrate."[10] This is a representative description for the complexity of prime distribution. But according to the famous "Langlands program", a difficult problem in one area could always be converted into an easy problem in other areas. Many recent efforts for this topic are all cross research of number theory and dynamic system theory, statistical theory, ergodic theory, etc. We believe the knowledge about the detailed distribution of primes will appear in elementary book of number theory soon.

# Reference

[1] M. Wolf, Physica 160A (1989) 24.

[2] Wang Liang, Huang Yan, Dai Zhi-cheng, Fractal in the statistics of Goldbach partition, http://arxiv.org/abs/nlin.CD/0601024

[3] Z. Gamba, J. Hernando, L. Romanelli, Phys. Lett. A 145 (1990) 106.

[4]M. Wolf, 1/f noise in the distribution of prime numbers, Physica A 241 (1997)

[5] S.R. Dahmen, S.D. Prado, T. Stuermer-Daitx, Similarity in the statistics of prime numbers, Physica A 296 (2001) 523–528.

[6]George G. Szpiro, The gaps between the gaps: some patterns in the prime number sequence, Physica A 341 (2004) 607 – 617

[7] B. Kra, The Green-Tao Theorem on arithmetic progressions in the primes: an ergodic point of view, Bull. Amer. Math. Soc. 43 (2006), 3-23

[8] B. Green and T. Tao, The primes contain arbitrarily long arithmetic progressions, http://arxiv.org/abs/math.NT/0404188

[9] B.-L. Hao, Elementary Symbolic Dynamics and Chaos in Dissipative Systems (World Scientific, Singapore, 1989).

[10]Leonard Euler, in G. Simmons, Calculus Gems, McGraw-Hill, New York, 1992